\documentclass[11pt]{amsart}
\usepackage{amssymb, paralist, xspace, graphicx, url, amscd, euscript, mathrsfs,stmaryrd,epic,eepic,color}
\usepackage[all]{xy}
\SelectTips{cm}{}

\numberwithin{equation}{section}

1

\numberwithin{subsection}{section}

\allowdisplaybreaks[1]

\newtheorem*{namedtheorem}{\theoremname}
\newcommand{\theoremname}{testing}

\newtheorem{theorem}[subsection]{Theorem}
\newtheorem{proposition}[subsection]{Proposition}
\newtheorem{proposition-definition}[subsection]
{Proposition-Definition}
\newtheorem{corollary}[subsection]{Corollary}
\newtheorem{lemma}[subsection]{Lemma}

\theoremstyle{definition}
\newtheorem{definition}[subsection]{Definition}
\newtheorem{notation}[subsection]{Notation}

\newtheorem{remark}[subsection]{Remark}

\newtheorem{problem}[subsection]{Problem}

\theoremstyle{remark}

\newcommand\cB{\mathcal{B}}
\newcommand\cC{\mathcal{C}}

\newcommand\cK{\mathcal{K}}
\newcommand\cL{\mathcal{L}}
\newcommand\cM{\mathcal{M}}

\newcommand\cO{\mathcal{O}}

\newcommand\ocM{\overline{\mathcal{M}}}
\newcommand\ocN{\overline{\mathcal{N}}}

\newcommand\ucK{\underline{\mathcal{K}}}

\newcommand\uB{\underline{B}}
\newcommand\uC{\underline{C}}
\newcommand\uD{\underline{D}}

\newcommand\uf{\underline{f}}

\newcommand\uG{\underline{G}}

\newcommand\uK{\underline{K}}

\newcommand\uS{\underline{S}}
\newcommand\uT{\underline{T}}

\newcommand\uW{\underline{W}}
\newcommand\uX{\underline{X}}
\newcommand\uY{\underline{Y}}
\newcommand\uZ{\underline{Z}}

\newcommand\NN{\mathbb{N}}

\newcommand\ZZ{\mathbb{Z}}

\newcommand\bfc{\mathbf{c}}

\newcommand\fB{\mathfrak{B}}
\newcommand\fC{\mathfrak{C}}

\newcommand\fM{\mathfrak{M}}

\newcommand\arr{\ifinner\to\else\longrightarrow\fi}

\def\displaytimes_#1{\mathrel{\mathop{\times}\limits_{#1}}}

\def\displayotimes_#1{\mathrel{\mathop{\bigotimes}\limits_{#1}}}

\newcommand\Aut{\operatorname{Aut}}

\newcommand\virt{{\operatorname{virt}}}

\newdir{ >}{{}*!/-5pt/@{>}}

\newcommand\double{\rightrightarrows}

\newcommand\doublelong[2]{\mathbin{\xymatrix{{}\ar@<3pt>[r]^{#1}
\ar@<-3pt>[r]_{#2}&}}}

\newlength{\ignora}

\newcommand{\B}{\mathop{\cB}\nolimits}

\theoremstyle{plain}

\newtheorem{lem}[subsection]{Lemma}
\newtheorem{prop}[subsection]{Proposition}

\theoremstyle{definition}

\newtheorem{rem}[subsection]{Remark}
\newtheorem{defn}[subsection]{Definition}
\newcommand{\Z}{\mathbb{Z}}

\numberwithin{equation}{subsection}

\begin{document}

\title{Stable logarithmic maps to Deligne--Faltings pairs II}

\author{Dan Abramovich}

\author{Qile Chen}

\address{Department of Mathematics\\
Brown University\\
Box 1917\\
Providence, RI 02912\\
U.S.A.}
\email{abrmovic@math.brown.edu, q.chen@math.brown.edu}


\thanks{Abramovich supported in part by NSF grants  DMS-0603284 and
  DMS-0901278. Chen  supported in part by funds from
  DMS-0901278.}
\begin{abstract}We make an observation which enables one to deduce the existence of an algebraic stack of log maps for all generalized Deligne--Faltings log structures (in particular simple normal crossings divisor) from the simplest case with characteristic generated by $\NN$ (essentially the smooth divisor case).
\end{abstract} 
\maketitle

\tableofcontents

\section{Introduction}
The idea of stable logarithmic maps was introduced in a legendary lecture by Bernd Siebert in 2001 \cite{Siebert}. However, the program has been on hold for a while, since Mark Gross and Bernd Siebert were working on other projects in mirror symmetry. Only recently they have taken up the unfinished project of Siebert jointly \cite{GS2}. The central object is a stack $\cK_\Gamma(Y)$ parameterizing what one calls {\em stable logarithmic maps} of log-smooth curves into a logarithmic scheme $Y$ with $\Gamma$ indicating the relevant numerical data, such as genus, marked points, curve class and other indicators related to the logarithmic structure. One needs to show $\cK_\Gamma(Y)$ is algebraic and proper. Gross and Siebert's approach builds on insights from tropical geometry, obtained by probing the stack of log maps using the standard log point. It covers the case of fixed targets with Zariski log structures, and is currently being extended to targets with relatively coherent log structures.

In \cite{Chen}, the second author considers another combinatorial construction of the stack $\cK_\Gamma(Y)$ when the logarithmic structure $Y$ on the underlying scheme $\uY$ is associated to the choice of a line bundle with a section. The motivating case is that of a pair $(\uY, \uD)$, where $\uD$ is a smooth divisor in the smooth locus of the scheme $\uY$ underlying $Y$. This particular situation enables him to approach the degeneration formula of \cite{Li-Ruan, Ionel-Parker, Li-degeneration}  in terms of logarithmic structures. It should be pointed out that these stable Logarithmic maps are not identical to those of Kim \cite{Kim}, though they are closely related.

Our point is that based solely on this special case, one can give a ``pure thought" proof of algebraicity and properness of the stack $\cK_\Gamma(Y)$ whenever $Y$ is a so called {\em generalized Deligne-Faltings} logarithmic structure. By saying $Y$ is a generalized Deligne-Faltings log structure we mean that there is a fine saturated sharp monoid $P$ and a sheaf homomorphism $P \to \ocM_{Y}$ which locally lifts to a chart  $P \to \cM_{Y}$; the slightly simpler {\em Deligne-Faltings log structure} is the case where $P = \NN^k$. This in turn covers many of the cases of interest, such as a variety with a simple normal crossings divisor, or a  simple normal crossings degeneration of a variety with a simple normal crossings divisor. We generalize it a bit further in Proposition \ref{Prop:moreDF}. It does not cover the case of a normal crossings divisor which is not simple, but we expect to cover this case using descent arguments. 

The purpose of this note is to set up a general categorical framework which enables us to make this construction. This general setup is of use not only for  $\cK_\Gamma(Y)$. In particular we have applications, pursued elsewhere \cite{ACGM}, to constructing the target of  {\em evaluation maps} of logarithmic Gromov-Witten theory.

All logarithmic schemes in this note are assumed to be fine and saturated logarithmic schemes - abbreviated fs log schemes -  unless indicated otherwise.

\subsection{Acknowledgements} We thank Danny Gillam, Mark Gross, Davesh Maulik,
Martin Olsson, Bernd Siebert, and Angelo Vistoli for numerous helpful comments
on this work.

\section{Logarithmic maps: a tale of two categories}\label{sec:log-map}
\subsection{Stable maps}
Let  $\uY$ be a projective scheme. The stack of stable maps to $\uY$ is defined as follows: one fixes discrete data $\Gamma=(g,n,\beta)$ where $g,n$ are non-negative integers standing for genus and number of marked points, and $\beta$ is a curve class on $\uY$. A {\em pre-stable} map to $\uY$ over a scheme $\uS$ is a diagram 
$$\xymatrix{\uC \ar[r]\ar[d]& \uY\\ \uS}
$$ 
where $\uC\to \uS$ is a proper flat family of n-pointed prestable curves, and $\uC \to \uY$ a morphism. The prestable map is {\em stable} if  on the fibers the groups $\Aut_{\uY}(\uC_s)$ are finite. Morphisms of prestable curves are defined as cartesian diagrams.

One easily sees that prestable maps form a category fibered in groupoids over the category of schemes. It is an important theorem that this fibered category is represented by an algebraic stack $\uK^{pre}(\uY)$, and the substack $\uK_\Gamma(\uY)$ of stable maps of type $\Gamma$ is proper with projective coarse moduli space \cite{Kontsevich}. When $\uY$ is smooth, there is a perfect obstruction theory, giving rise to a virtual fundamental class $[\uK_\Gamma(\uY)]^\virt$ underlying the usual algebraic treatment of Gromov--Witten theory \cite{Li-Tian, Behrend-Fantechi}.

The main result of this note is an analogue of the following evident result: assume $\uY = \uY_1\times_{\uY_2} \uY_3$. Then $$\uK^{pre}(\uY) = \uK^{pre}(\uY_1)\mathop{\times}\limits_{\uK^{pre}(\uY_2)}\uK^{pre}(\uY_3).$$ 

\subsection{Stable logarithmic maps as a stack over $\mathfrak{LogSch}$.}\label{ss:map-stack-over-logsch}
Let $Y$ be a  fs logarithmic scheme with projective underlying scheme $\uY$.
One can repeat the construction above, replacing prestable curves by {\em log smooth} curves \cite{FKato}, and replacing all morphisms of schemes by morphisms of log schemes: a {\em pre-stable log map over $S$} is a diagram of {\em logarithmic} schemes 
$$\xymatrix{C \ar[r]\ar[d]& Y\\ S}
$$ 
where $C \to S$ is a log smooth curve and $C \to Y$ a morphism of log schemes. We define such a map to be {\em stable} if the underlying prestable map is stable. Arrows are defined using cartesian diagrams: 
$$\xymatrix{C'\ar[r]\ar[d]&C \ar[r]\ar[d]& Y\\S' \ar[r]& S}
$$ 

Again it is evidently a category fibered in groupoids, but this time {\em over the category $\mathfrak{LogSch}$ of fs logarithmic schemes}. It is proven in \cite{Chen} when the log structure $Y$ is given by a line bundle with a section, and more general in \cite{GS2}, that this category is an fs logarithmic algebraic stack: there is a logarithmic algebraic  stack $\cK_\Gamma(Y) = (\ucK,\cM_{{\ucK}})$, where the log structure is fs, such that stable logarithmic maps over $S$ are equivalent to log morphisms $S \to \cK_\Gamma(Y)$. Denote the universal log smooth curve by $\fC\to \cK_\Gamma(Y)$.

It is natural to search for general criteria for algebraicity of such logarithmic moduli which is analogous to Artin's work \cite{Artin}. We do not address this general question here.

\subsection{Stable logarithmic maps as a stack over $\mathfrak{Sch}$.} The existence of $\cK_\Gamma(Y)$ has immediate strong implications on the structure of stable log maps.  Objects of the underlying  stack $\ucK_\Gamma(Y)$ over a scheme $\uS$ can be understood as follows: an object is after all an arrow $\uS \to \ucK_\Gamma(Y)$. It automatically gives rise to 
a cartesian diagram
$$\xymatrix{C^{\min} \ar[r]\ar[d]& \fC \ar[r]\ar[d]                &Y\\ 
                     S^{\min}  \ar[r]\ar[d]&\cK_\Gamma(Y) \ar[d]\\
                     \uS\ar[r]        & \ucK_\Gamma(Y),}
$$ 
in particular an object $S ^{\min} \to \cK_\Gamma(Y)$, but here {\em the logarithmic structure  $S ^{\min}$ is pulled back from   $\cK_\Gamma(Y)$}. Moreover, {\em every} stable log map factors uniquely through one of this type: Given a stable log map over $S$ we have a morphism $S \to \cK_\Gamma(Y)$ by definition, giving rise to an extended  diagram
$$\xymatrix{C \ar[r]\ar[d] &C^{\min} \ar[r]\ar[d]& \fC \ar[r]\ar[d]                &Y\\ 
                     S \ar[r] &S^{\min}  \ar[r]\ar[d]&\cK_\Gamma(Y) \ar[d]\\
                     &\uS\ar[r]        & \ucK_\Gamma(Y).}
$$ 
Following B. Kim \cite{Kim} we call a log map over $S$ {\em minimal} (not to be confused with log minimal models of the minimal model program) if $S \to \cK_\Gamma(Y)$ is strict, namely the log structure on $S$ is the pullback of the log structure on $\cK_\Gamma(Y)$. It follows tautologically that the underlying stack $\ucK_\Gamma(Y)$ precisely parametrizes log maps with minimal log structure.

In fact this thought process is reversible: the construction of \cite{Chen} in the case of a Deligne-Faltings log structure of rank 1 goes by way of constructing a proposed minimal log structure associated to any log map, and verifying that log maps where the log structure is the proposed minimal one are indeed minimal (every object  maps uniquely to a minimal one), and form an algebraic stack over $\mathfrak{Sch}$ carrying a logarithmic structure. 
In short, the second categorical interpretation, of $\ucK_\Gamma(Y)$ as a stack over $\mathfrak{Sch}$, takes precedence here.

One is tempted to try to mimic the same construction in general. {\em This is not the route taken here.} In fact we use the universality of the category $\cK_\Gamma(Y)$ over $\mathfrak{LogSch}$ for given $Y$ to deduce its algebraicity from cases of simpler $Y$. Minimal object are obtained as an afterthought in Proposition \ref{prop:univ-min-map} and described combinatorially in Section \ref{ss:DF-marked-graph}. It is worthwhile setting this up in general.

\subsection{The general setup} Consider a commutative diagram 

$$
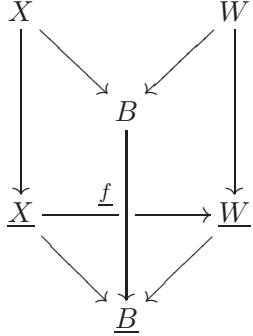
\xymatrix{
X \ar[dd]\ar[dr] && W\ar[dd]\ar[dl]\\
& B\ar[dd]\\
\uX\ar@{->}'[r]^(.8)\uf[rr]\ar[dr]  && \uW\ar[dl] \\
&\uB
}$$
where $X \to B$ is an integral morphism of fs log schemes  with $\uX \to \uB$ flat and proper,  $W$ is an fs log scheme, and $\uX \to \uW$ a morphism. We can define a category $\text{Lift}_{\uf}$ fibered in groupoids over $\mathfrak{LogSch}$ whose objects over $S$ are morphisms $S \to B$ of fs log schemes  with a lifting $f_S:X_S \to W_S$ of the underlying morphism  $\uf_S:\uX_{\uS} \to \uW_{\uS}$. The arrows are again defined by taking cartesian diagrams. We can ask the following question:

\begin{problem}\label{problem}
Is the category $\text{Lift}_{\uf}$ equivalent to an fs logarithmic algebraic stack $\cL ift_{\uf}$? Under what conditions is it proper?
\end{problem}

In fact if convenient we can remove the geometry of $W$ entirely from the picture, mimicking the methods of Olsson \cite{Olsson}:  consider $\overline{W} = \uW \times_{\uB} B$ and  $Z = X \times_{\overline W} W$. Then $Z$ is a logarithmic scheme over $X$. The category $\text{Lift}_{\uf}$ is evidently equivalent to the category $\text{Sec}_{X/B}(Z/X)$ whose objects over a log morphism $S \to B$ are sections $X_S \to Z_S$ of $Z_S\to X_S$. 

The main result is the following: 
\begin{theorem}\label{Th:Sec-limits}
Let $\Delta = (Z_\alpha, \pi_{\alpha\beta}: Z_{\alpha}\to Z_{\beta})$ be a finite diagram of fs log schemes over $X$, with final object $X$.  Assume $$Z = \varprojlim(\Delta)$$ in the category of fs log schemes. Then 
$$\text{Sec}_{X/B}(Z/X) = \varprojlim \text{Sec}_{X/B}(Z_\alpha/X)
$$ namely it is the limit in the category of fs log schemes.
\end{theorem}

If the reader finds general categorical limits a bit off-putting, the main cases needed for our applications are (1) fiber products, and (2) equalizers, which can be described as fiber products using the diagonal morphisms.


The existence of such limits in the category of fs log schemes is proven  in \cite{KKato}: the case of arbitrary log structures is treated in (1.6), coherent log structures in (2.6), and fine log structures follow from (2.7); the case of fs log structures follows from the analogous adjoint functor with $P^{int}$ replaced by $P^{sat}$ (\cite[Chapter I, Proposition 1.2.3]{Ogus}).

\begin{proof}[Proof of the theorem]

An object of $\text{Sec}_{X/B}(Z/X)$ over an arrow $S\to B$ is by definition a section $s:X_S \to Z$, and composing with the canonical maps we get $s_\alpha:X_S \to Z_\alpha$ such that for each arrow $\pi_{\alpha\beta}: Z_\alpha \to Z_\beta$ in $\Delta$ we have $\pi_{\alpha\beta}\circ s_{\alpha} = s_{\beta}$. This in particular gives us a diagram of objects $s_\alpha\in\text{Sec}_{X/B}(Z_\alpha/X)(S)$ with $\pi_{\alpha\beta}(s_{\alpha}) = s_{\beta}$, namely an object of $\varprojlim\text{Sec}_{X/B}(Z_i/X)$. The correspondence on the level of arrows is similar (though maybe more confusing).

The process is completely reversible, hence the equivalence.
\end{proof}

\section{The stacks of stable logarithmic maps}

Theorem \ref{Th:Sec-limits} applies directly to the category $\text{Lift}_{\uf}$ and Problem \ref{problem}. In the present section we make this as explicit as possible, including a discussion of {\em contact orders}, the deformation-invariant numerical data encoded in the logarithmic structure.

\subsection{The case of Deligne-Faltings pairs: setup}
Let the target $Y$ be a DF pair. We first break up the log structure $\cM_{Y}$ into the known cases as in \cite{Chen}. We may assume that $P\cong \NN^{m}$, and index the $m$ copies of $\NN$ by the set $\{1,2,\cdots,m\}$. Note that the map $P\to\ocM_{Y}$ locally lifts to a chart. Consider the $i$-th copy of $\NN$, and denote by $\NN_{i}\hookrightarrow P$. Locally the following composition
\[\NN_{i}\to P\to \cM_{Y}\]
defines a rank 1 DF-log structure $\cM_{i}$ on $\underline{Y}$. Note that $\cM_{i}\subset \cM_{Y}$ is a sub-log structure, and the following decomposition holds:
\begin{equation}\label{equ:decomp-DF-log}
\cM_{Y}=\cM_{1}\oplus_{\cO_{Y}^{*}}\cM_{2}\oplus_{\cO_{Y}^{*}}\cdots\oplus_{\cO_{Y}^{*}}\cM_{m}.
\end{equation}

Denote by $Y_{i}=(\underline{Y},\cM_{i})$ for all $i$. We view $\underline{Y}$ as a log scheme with trivial log structure. The above decomposition is equivalent to 
\begin{equation}\label{equ:decomp-DF-pair}
Y = Y_{1}\times_{\underline{Y}}Y_{2}\times_{\underline{Y}}\cdots\times_{\underline{Y}}Y_{m},
\end{equation}
which can be written as 
\[Y = \varprojlim Y_{i}.\]
Note that each $\cM_{i}$ corresponds to a pair $(L_{i},s_{i})$, consisting of a line bundle $L_{i}$ with a section $s_{i}\in H^{0}(L_{i}^{\vee})$.

\subsection{Numerical data and contact orders}
We introduce a notation for the numerical data:

\begin{notation}\label{nota:DF-dis-data}
Denote by $\Gamma=(\beta,g,n,\bfc)$ where
\begin{enumerate}
 \item $\beta\in H^{2}(\uY,\ZZ)$ is a curve class;
 \item $g$ is a non-negative integer, which will denote the genus of the source curve;
 \item $n$ is a non-negative integer, which will denote the number of marked points;
 \item $\bfc=\{(c_{i}^{j})_{i=1}^{m}\}_{j=1}^{n}$ is a set of tuples, where $c_{i}^{j}$ denotes the contact order of the $j$-th marking with respect to $\cM_{i}$. Hence we require that 
 \[\beta\cdot c_{1}(L^{\vee}_{i})=\sum_{j=1}^{n}c_{i}^{j},\]
 where $c_{1}(L^{\vee}_{i})$ is the first chern class of $L^{\vee}_{i}$.
\end{enumerate}
\end{notation}

Note that the set $\bfc$, hence the stack $\cK_{\Gamma}(X)$ we will construct later, depends on the choices of $P$. For example, consider the union of two disjoint divisors $D=D_{1}\cup D_{2}$ in $\uY$. Denote by $I_{i}$ the ideal sheaf of $D_{i}$, and $s_{i}:I_{i}\to \cO_{\uY}$ the natural injection. Let $\cM_{Y}$ be the standard log structure associated to $D$ on $\uY$, and $Y$ the log scheme $(\uY,\cM_{Y})$. Now we have two map of sheaves of monoids 
\begin{equation}\label{equ:diff-present}
\NN\to \ocM_{Y}, \mbox{\ \ \ and \ \ \ }\NN^{2}\to \ocM_{Y},
\end{equation}
where the first one is given by the ideal sheaf $I_{D}$ of $D$ with the natural injection $s_{D}:I_{D}\to \cO_{Y}$, and the second one is given by the two pairs $(I_{i},s_{i})$ for $i=1,2$. Note that both maps of monoids in (\ref{equ:diff-present}) can be lifted to  charts of $\cM_{Y}$ locally. 

Consider a marked point $p$ with assigned contact order $c$ under the choice $P=\NN$. This includes the following  two cases: $p$ can have contact orders $c$ with $D_{1}$ and $0$ with $D_{2}$, or $c$ with $D_{2}$ and $0$ with $D_{1}$. However if we choose the second case where $P = \NN^{2}$, then the contact orders along both $D_{1}$ and $D_{2}$ will be specified. Thus, if we choose different choice of $P$, hence different $\bfc$, we will obtain different tangency conditions, hence different stacks $\cK_{\Gamma}(Y)$.

If we consider the stack $\cK_{g,n}^{log}(Y,\beta)$ as in \cite{Chen}, which parametrizes minimal stable log maps with all possible choices of contact orders. This does not depend on the choice of the monoid $P$. For any choice of $P$ and $\Gamma$, the stack $\cK_{\Gamma}(Y)$ is an open and closed substack of $\cK_{g,n}^{log}(Y,\beta)$. 


A further refined and entirely canonical formalism of contact orders follows from \cite{ACGM}: given any fs logarithmic scheme $Y$ one defines an Artin stack $\wedge Y$,  locally of finite type over $\uY$, parameterizing standard log points in $Y$. We call it  {\em the evaluation space of $Y$.} Given a stable log map $f:C \to Y$ and an integer $j$ with $1\leq j\leq n$, the restriction $f_{\Sigma_j}: \Sigma_j\to Y$ of $f$ to the $j$-th marking, is an element of $\wedge Y$. This defines the $j$-th evaluation map $\cK_{g,n}(Y,\beta) \to \wedge Y$.

\begin{definition}\label{nota:can-contact}
A {\em logarithmic sector} of $Y$ is an element $c\in \pi_0(\wedge Y)$, namely a connected component of  the evaluation space of $Y$. A stable map $f$ is said to have canonical contact order $\bfc=(c^1,\ldots,c^n)$ if the  $j$-th evaluation map lands in the logarithmic sector $c^j\in \pi_0(\wedge Y)$. 
\end{definition} 

In \cite[Formula (1.6)]{GS2}, Gross and Siebert introduce numerical
  data denoted $u_p$. This data is closely related to our canonical
  contact orders; this  relationship is made more precise in \cite{ACGM}.




The discussion which follows works equally well with the explicit and down-to-earth  contact orders of Definition \ref{nota:DF-dis-data}, as with the canonical contact orders of Definition \ref{nota:can-contact}.

\subsection{The case of Deligne-Faltings pairs: the stacks}

\begin{defn}
Denote by $\cK_{\Gamma}(Y)$ (respectively $\cK_{\Gamma}^{pre}(Y)$) the category fibered over $\mathfrak{LogSch}$, which parametrizes stable Log maps (respectively pre-stable log maps) with discrete data $\Gamma$ as in Notation \ref{nota:DF-dis-data}.
\end{defn}

Note that for each $i$, we have a set $\bfc_{i}=\{c_{i}^{j}\}_{j=1}^{n}$. Consider the discrete data $\Gamma_{i}=(\beta,g,n,\bfc_{i})$. Then by \cite{Chen} we have the log stack $\cK_{\Gamma_{i}}(Y_{i})$, parameterizing minimal stable Log maps with discrete data $\Gamma_{i}$. (If instead we use canonical contact orders, then $\bfc_i$ is simply the image of $\bfc$ under the canonical map $\pi_0(\wedge Y) \to \pi_0(\wedge Y_i)$.)

Consider the stack $\cK_{g,n}(\uY,\beta)$ of usual stable maps. Since the universal curve is prestable, it has a canonical log smooth structure coming from $\fM_{g,n}$ using F. Kato's theorem \cite{FKato}. We put this in the setup of $\text{Sec}_{X/B}(Z/X)$  by setting $B$ to be $\cK_{g,n}(\uY,\beta)$ or $\cK^{pre}_{g,n}(\uY,\beta)$ with the pull-back log structure coming from $\fM_{g,n}$; $X$ to be the universal curve with its canonical log structure; and $W=Y\times B$. Consider $Z = X \times_{(\uW\times_{\uB} B)} W$ as before. Denote by $Z_{i} = X \times_{(\uW\times_{\uB}B)} W_{i}$ for $i=1,\cdots,m$. Let $Z_{0}=X \times_{(\uW\times_{\uB} B)} \uW$. Then by (\ref{equ:decomp-DF-pair}), we have 
\[Z\cong Z_{1}\times_{Z_{0}}Z_{2}\times_{\Z_{0}}\cdots\times_{\Z_{0}}Z_{m},\]
or equivalently
\[Z = \varprojlim Z_{i}.\]

\begin{corollary}\label{cor:target-DF}
The fibered categories $\cK_{\Gamma}(Y)$ and $\cK^{pre}_{\Gamma}(Y)$ are represented by algebraic stacks equipped with fs log structures. Furthermore, $\cK_{\Gamma}(Y)$ is representable and finite over $\cK_{g,n}(\uY,\beta)$.
\end{corollary}
\begin{proof}
For the first statement, it follows from Theorem \ref{Th:Sec-limits} that 
\begin{equation}
\cK_{\Gamma}^{pre}(Y)\cong \cK_{\Gamma_{1}}^{pre}(Y_{1})\times_{\cK_{g,n}^{pre}(\uY,\beta)}\cdots\times_{\cK_{g,n}^{pre}(\uY,\beta)}\cK_{\Gamma_{m}}^{pre}(Y_{m}),
\end{equation}
where the fibered product are taken in the category of fs log schemes. It was shown in \cite{Chen} that the first statement holds for $\cK_{\Gamma_{i}}^{pre}(Y_{i})$. Note that 
\[\cK_{\Gamma}(Y)\cong \cK_{\Gamma}^{pre}(Y) \times_{\cK_{g,n}^{pre}(\uY,\beta)}\cK_{\Gamma}(\uY,\beta).\] 
Thus, the first statement follows. 

Similarly, we have a fiber product of fs log stacks
\begin{equation}\label{equ:DF-stack}
\cK_{\Gamma}(Y)\cong \cK_{\Gamma_{1}}(Y_{1})\times_{\cK_{g,n}(\uY,\beta)}\cdots\times_{\cK_{g,n}(\uY,\beta)}\cK_{\Gamma_{m}}(Y_{m}).
\end{equation}
Now the second statement follows from the finiteness of taking product in the category of fs log schemes, see \cite[Chapter 2, 2.4.5]{Ogus}. 
\end{proof}


\subsection{The case of generalized Deligne-Faltings pairs}
Slightly generalizing the above situation, we first give the following:

\begin{defn}\label{defn:gen-DF-log}
A log structure $\cM_{Y}$ on $\uY$ is called a {\em generalized Deligne-Faltings log structure}, if there exists a toric sharp monoid $P$ and a map $P\to \ocM_{Y}$ from the sheaf of constant monoid $P$, which locally lifts to a chart. The log scheme $Y=(\uY,\cM_{Y})$ is called a {\em generalized Deligne-Faltings (DF) pair}.
\end{defn}

Consider the target $Y$ with generalized DF log structure $\cM_{Y}$, with a fixed map $P\to \ocM_{Y}$ as in Definition \ref{defn:gen-DF-log}. Again we have $B = \cK_{g,n}(\uY,\beta)$ with the pull-back log structure coming from $\fM_{g,n}$; $X$ is the universal curve of $B$ with its canonical log structure; $W=Y\times B$; and $Z = X \times_{(\uW\times_{\uB}B)} W$ as before. 

By \cite[Chapter 1, 2.1.9(7)]{Ogus}, we have the following coequalizer diagram of monoids
\begin{equation}\label{equ:monoid-coeq}
\xymatrix{
\NN^{n_{2}} \ar@<.5ex>[r]^{v_{1}} \ar@<-.5ex>[r]_{v_{2}} & \NN^{n_{1}} \ar[r]^{q} & P,
}
\end{equation}
where $n_{1}$ and $n_{2}$ are non-negative integers. This is equivalent to say that
$$P = \lim_\rightarrow(\NN^{n_{2}} \double \NN^{n_{1}}).$$
Thus we have the following push-out diagram of fs monoids:
\begin{equation}\label{diag:monoid-pushout}
\xymatrix{
\NN^{n_{2}}\oplus\NN^{n_{2}} \ar[rr]^{id\oplus id} \ar[d]_{v_{1}\oplus v_{2}} && \NN^{n_{2}} \ar[d] \\
\NN^{n_{1}} \ar[rr] && P.
}
\end{equation}
 
Since the map $P\to \ocM_{Y}$ locally lifts to a chart, the two compositions
\[ \NN^{n_{1}}\to P \to \cM_{Y} \mbox{\ \ and\ \ } \NN^{n_{2}}\to P \to \cM_{Y}\]
induce two DF log structures $\cM_{1}$ and $\cM_{2}$ respectively, with the following coequalizer diagram of log structures:
\begin{equation}\label{diag:coequal-log}
\xymatrix{
\cM_{2} \ar@<0.5ex>[r]^{v_{1}} \ar@<-.5ex>[r]_{v_{2}} & \cM_{1} \ar[r]^{q} & \cM_{Y}.
}
\end{equation}
Denote by $\cM_{3}= \cM_{2}\oplus_{\cO^{*}_{Y}} \cM_{2}$. Since $\cM_{2}$ is DF, there is a natural map 
\[\NN^{n_{2}}\oplus\NN^{n_{2}} \to \ocM_{3}\]
which locally lifts to a chart. Thus $\cM_{3}$ is also a DF log structure. It is not hard to check that the coequalizer diagram (\ref{diag:coequal-log}) is equivalent to the following push-out diagram:
\begin{equation}\label{diag:limit-pushout}
\xymatrix{
\cM_{3} \ar[rr]^{id\oplus id} \ar[d]_{v_{1}\oplus v_{2}} && \cM_{2} \ar[d] \\
\cM_{1} \ar[rr] && \cM_{Y}.
}
\end{equation}  

Denote by $Y_{i}=(\underline{Y},\cM_{i})$, for $i=1,2,$ and $3$. Then (\ref{diag:coequal-log}) induces a equalizer of log schemes:
\[
Y\to Y_{1} \double Y_{2},
\]
or equivalently
\begin{equation}\label{equ:log-sch-lim}
Y= \varprojlim Y_i.
\end{equation}
By construction, the map $\NN^{n_{i}}\to \ocM_{i}$ locally lifts to a chart. Hence each $Y_{i}$ is a DF pair.

The diagram (\ref{diag:limit-pushout}) implies that (\ref{equ:log-sch-lim}) is equivalent to the following cartesian diagram of fs log schemes:
\begin{equation}\label{diag:log-sch-car}
\xymatrix{
Y \ar[r]^{} \ar[d]_{} & Y_{2} \ar[d]^{}\\
Y_{1} \ar[r]_{} & Y_{3}.
}
\end{equation}
Denote by $W_{i}= Y_{i}\times B$, and $Z_{i}=X \times_{(\uW_{i}\times_{\uB}B)} W_{i}$. Then (\ref{equ:log-sch-lim}) and (\ref{diag:log-sch-car}) implies that 
\begin{equation}\label{equ:log-sch-lim-Z}
Z= \varprojlim Z_i,
\end{equation}
or equivalently the cartesian diagram of fs log schemes:
\begin{equation}\label{diag:log-sch-car-Z}
\xymatrix{
Z \ar[r]^{} \ar[d]_{} & Z_{2} \ar[d]^{}\\
Z_{1} \ar[r]_{} & Z_{3}.
}
\end{equation}

We can define discrete data using canonical contact orders $\bfc$ as in Definition \ref{nota:can-contact}. Alternatively, 
we can define a concret and down-to-earth version  using $Irr(P)$, the set of irreducible elements in $P$. Since the monoid $P$ is toric and sharp, $Irr(P)$ forms a set of generators of $P$. For each element $\alpha\in Irr(P)$, there is a free sub-monoid $\NN_{\alpha}\cong \NN \hookrightarrow P$, which induces a rank 1 DF sub-log structure $\cM_{\alpha}\subset \cM_{Y}$. Denote by $(L_{\alpha},s_{\alpha})$ to be the pair of line bundle and sections corresponding to $\cM_{\alpha}$. Again, we fix the numerical data $\Gamma=(\beta,g,n,\bfc)$ such that 
\begin{enumerate}
 \item the data $\beta, g,n $ denote the curve class, genus, and number of marked points as in Notation \ref{nota:DF-dis-data};
 \item $\bfc=\{(c_{\alpha}^{j})_{\alpha\in Irr(P)}\}_{j=1}^{n}$ is a set of tuples, where $c_{\alpha}^{j}$ denotes the contact order of the $j$-th marking with respect to $\cM_{\alpha}$ satisfying 
 \[\beta\cdot c_{1}(L^{\vee}_{\alpha})=\sum_{j=1}^{n}c_{\alpha}^{j}.\]
\end{enumerate} 
 
Note that the data $\Gamma$ induces discrete data $\Gamma_{i}$ for log maps to $Y_{i}$ via (\ref{diag:monoid-pushout}). We obtain the following:

\begin{corollary}\label{Cor:stable-maps}
Assume that the logarithmic structure of $Y$ is a generalized Deligne--Faltings log structure. Then the categories $\cK_{\Gamma}(Y)$ and $\cK_{\Gamma}^{pre}(Y)$ are algebraic stacks with fs log structures. Furthermore, the stack $\cK_{\Gamma}(Y)$ is representable and finite over $\cK_{g,n}(\uY,\beta)$.
\end{corollary}
\begin{proof}
Again Theorem \ref{Th:Sec-limits} implies that 
\begin{equation}
\cK_{\Gamma}^{pre}(Y)\cong \cK^{pre}_{\Gamma_{1}}(Y_{1})\times_{\cK^{pre}_{\Gamma_{3}}(Y_{3})}\cK^{pre}_{\Gamma_{2}}(Y_{2})
\end{equation}
and 
\begin{equation}\label{equ:gen-DF-stack}
\cK_{\Gamma}(Y)\cong \cK_{\Gamma_{1}}(Y_{1})\times_{\cK_{\Gamma_{3}}(Y_{3})}\cK_{\Gamma_{2}}(Y_{2})
\end{equation}
where the fibered product is taken over the category of fs log schemes. Now the statement follows from the same argument as in Corollary \ref{cor:target-DF}. 
\end{proof} 

\begin{remark}\label{rem:all-contact}
Denote by $\cK_{g,n}(Y,\beta)=\coprod_{\Gamma\in \Lambda}\cK_{\Gamma}(Y)$, where $\Lambda$ is the set of all possible contact orders $\Gamma=(\beta,g,n,\bfc)$ with fixed $g,n$ and $\beta$. Then $\cK_{g,n}(Y,\beta)$ represents the category of stable log maps over $\mathfrak{LogSch}$ with all possible contact orders. Similarly, we can define the category $\cK_{g,n}^{pre}(Y,\beta)=\coprod_{\Gamma\in \Lambda}\cK_{\Gamma}^{pre}(Y)$, which parametrizes all pre-stable log maps with all possible contact orders. Clearly, this is also represented by an algebraic stack with natural fs log structure.
\end{remark}

\subsection{A further  generalization}

We can weaken the Deligne--Faltings assumption in Corollary \ref{Cor:stable-maps} as follows. Consider a fine and saturated log scheme $Y$ with a surjective homomorphism of sheaves of monoids $\NN_{Y}^{k}\to \ocM_{Y}$. This induces a new log scheme $Y'$ with underlying scheme $\uY$, and the log structure $\cM_{Y'}$, which is locally given by the log structure associated to the composition of $\NN_{Y}^{k}\to \ocM_{Y}$ with a local chart $\ocM_{Y} \to \cM_{Y} \to \cO_{\uY}$. Thus we have a natural map $\psi: Y\to Y'$ with $\underline{\psi}=id_{\uY}$. Note that the log structure on $Y'$ is Deligne-Faltings.

As in Section \ref{ss:map-stack-over-logsch}, we have a fibered category $\cK_{g,n}(Y,\beta)$ parameterizing genus $g$, $n$-pointed stable log maps with curve class $\beta$ to $Y$ over $\mathfrak{LogSch}$. One can also introduce canonical contact orders as in  Definition \ref{nota:can-contact}. Note that the map $\psi$ induces a natural map of the fibered categories 
\begin{equation}\label{map:log-embed}
\phi: \cK_{g,n}(Y,\beta)\to \cK_{g,n}(Y',\beta).
\end{equation}

\begin{proposition}\label{Prop:moreDF}
The fibered category $\cK_{g,n}(Y,\beta)$ is a logarithmic algebraic stack. Furthermore, the underlying morphism $\underline{\phi}$ is representable and finite, and the map of log structures $\phi^*\cM_{\cK_{g,n}(Y',\beta)} \to \cM_{\cK_{g,n}(Y,\beta)}$ is surjective.
\end{proposition} 

\begin{proof}
By Corollary \ref{Cor:stable-maps}, we have that $\cK_{g,n}(Y',\beta)$ is a proper logarithmic algebraic stack. It therefore suffices to prove the second statement in the proposition.

We construct $\cK_{g,n}(Y,\beta)$ locally over $\cK_{g,n}(Y',\beta)$. We may assume we have a fs log scheme of finite type $B'$ and a stable log map $(C'/B', C' \to Y')$. We claim that $B :=B' \times_{\cK_\Gamma(Y')}\cK_\Gamma(Y) \to B'$ is a finite map of the underlying schemes, and the log structure of $B'$ surjects to that of $B$.

We take $Z= Y \times_{Y'} C'$ to be the fiber product in the category of fine log schemes. Then the canonical projection $Z\to C'$ induces a closed embedding $\underline Z \to \uC'$. By \cite[VIII Theorem 6.5]{SGA3} and the reduction argument \cite[Theorem 6(3)]{Abramovich}, there is a universal closed sub-scheme of $\uW\subset \uB'$ satisfying that for any $\uT \to \B'$ if $\uZ\times_{\uB'}\uT \to \uC'\times_{\uB'}\uT$ is an isomorphism, then the map $\uT \to \B'$ factor through $\uW$ uniquely. Note that any element of $B$ over some fs log scheme $T$ is a commutative diagram
\begin{equation}\label{diag:new-object}
\xymatrix{
Z_{T} \ar[r] \ar[d] & Z \ar[r] \ar[d] & Y \ar[d] \\
C'_{T} \ar[r] \ar[d] \ar@/^1pc/[u]^{s_{T}} & C' \ar[d] \ar[r] & Y' \\
T \ar[r] & B'.
}
\end{equation}
where the two upper squares are cartesian squares of fine log schemes, and the bottom one is a cartesian square of fs log schemes. The map $s_{T}$ is a section of $Z_{T}\to C'_{T}$. Thus we have an isomorphism of the underlying schemes $\uZ_{T}\cong \uC'_{T}$. This implies that the map $\uT\to \uB'$, hence $\uB\to \uB'$ factors through $\uW$. We replace $B'$ by $\uW$, with the pullback log structure.  

Since the problem is local, we may assume that there are sections $\sigma_i:  B' \to C'$ for $i=1,2,\cdots, n$ landing in the generic locus, and meeting every component of every fiber. Let $B'_{0}=B'$, $Z_{0}=Z$, and $C'_{0}=C'$. We construct $B$ by induction.

Denote by $B'_1= B'_{0}\times_{C'_{0}}Z_{0}$, where the product is taken via $\sigma_1$ in the category of fs log schemes. It follows that the underlying map of the first projection $h_{1}: B_{1}'\to B'_{0}$ is finite, and the map of log structures $h^{\flat}_{1}:\cM_{B'_{0}}\to \cM_{B_{1}'}$ is surjective. We claim that the map $B \to B'$ factors through $B'_1$ uniquely. To see this, we pick a commutative diagram as in (\ref{diag:new-object}). This induces a commutative diagram
\[
\xymatrix{
Z_{0,T} \ar[r] \ar[d] & Z_{0}  \ar[d] \\
C'_{0,T} \ar[r] \ar[d] \ar@/^1pc/[u]^{s_{T}} & C'_{0} \ar[d] \\
T \ar[r] \ar@/^1pc/[u]^{(\sigma_{i})_{T}} & B'_{0} \ar@/_1pc/[u]_{\sigma_{i}}.
}
\]
hence a commutative diagram
\[
\xymatrix{
T \ar@/^/[rrd] \ar@{-->}[rd] \ar@/_/[rdd] && \\
& B_{1}' \ar[r] \ar[d] & Z_{0} \ar[d]\\
& B'_{0}     \ar[r]^{\sigma_{1}} & C'_{0}
}
\]
where the dashed arrow is induced by the universal property of fibered product. This proves the claim.

Replacing $B'_{i-1}$ by $B'_{i}$ one at a time, and repeating the previous step, we obtain a sequence of maps $\{h_{i}:B_{i}' \to B_{i-1}'\}_{i}$ such that 
\begin{enumerate}
 \item the underlying map $\underline{h}_{i}$ is finite;
 \item the map of log structures $h_{i}^{\flat}:h_{i}^{*}\cM_{B_{i-1}'}\to \cM_{B_{i}'}$ is surjective;
 \item There is a canonical map $B\to B_{i-1}$ of fibered categories over $\mathfrak{LogSch}$ which factors through $B_{i}$ uniquely.
\end{enumerate}

Denote by $Z_{n}\to C'_{n} \to B'_{n}$ the pull-back of $Z\to C' \to B' $ via $B'_{n}\to B'$. By \cite[Lemma 3.5 iii]{Olsson-ENS}, after taking finitely many sections $\sigma_{i}$, we may assume $Z_{n}\to C'_{n} \to B'_{n}$ has the property that $Z_{n}\to C'_{n}$ is an isomorphism along the smooth non-marked locus of each fiber. 

We claim that $Z_{n} \to C'_{n}$ is an isomorphism of log schemes. Since the underlying map is an isomorphism, it is enough to show that the induced map of characteristics $\bar{l}:\ocM_{Z_{n}}\to \ocM_{C'_{n}}$ is an isomorphism. Note that we have a cartesian diagram of fine log schemes:
\[
\xymatrix{
Z_{n} \ar[rr] \ar[d] && Y \ar[d] \\
C'_{n}\ar[rr] && Y'.
}
\]
Since the map of log structures $\cM_{Y'}\to \cM_{Y}$ is surjective, the map of monoids $\bar{l}$ is also surjective. 

Denote by $\underline{f}: \uC'_{n}\to \uY$ the underlying stable map, and $\pi_{n}: C'_{n}\to B_{n}'$ the map of log curves. Locally at a marked point $p\in C'_{n}$, we have a push-out diagram of fine monoids:
\[
\xymatrix{
\underline{f}^{*}\ocM_{Y'} \ar[rr] \ar@{->>}[d] && \ocN_{p}\oplus\pi^{*}_{n}\ocM_{B_{n}'} \ar@{->>}[d]^{\bar{l}} \\
\underline{f}^{*}\ocM_{Y}  \ar[rr] && \ocM_{Z_{n}}
}
\]
where $\ocN_{p}\cong \NN =\ <\delta>$ is the standard log structure given by the marking associated to $p$. Pick two elements $c_{1} \delta + e_{1},\  c_{2} \delta + e_{2} \in \ocN_{p}\oplus\pi^{*}_{n}\ocM_{B_{n}'}$ locally at $p$, where $e_{1}, e_{2}\in \pi^{*}_{n}\ocM_{B_{n}'}$, and $c_{1},c_{2}$ are two non-negative integers. Assume that $\bar{l}(c_{1}\cdot \delta + e_{1})= \bar{l}(c_{2}\cdot \delta + e_{2})$. By generalizing to a nearby smooth non-marked point of $p$, we have $\bar{l}(e_{1})=\bar{l}(e_{2})$. Since $\bar{l}$ is an isomorphism at the generic smooth non-marked points, this implies that $e_{1}=e_{2}$. Thus, we have $\bar{l}(c_{1}\delta)=\bar{l}(c_{2}\delta)$. If $c_{1}\neq c_{2}$, then there exists another positive integer $c_{3}$ such that $\bar{l}(c_{3}\delta)=0$, which implies that $c_{3}\delta=0$ in $\ocN_{p}$. This is a contradiction. Therefore, we have $c_{1}\cdot \delta + e_{1} = c_{2}\cdot \delta + e_{2}$. This implies that $\bar{l}$ is an injection, hence an isomorphism at each marked point. A similar argument as above will imply that $\bar{l}$ is also an isomorphism at each node. Therefore, the map $Z_{n}\to C_{n}'$ is an isomorphism of log schemes.

Note that we have an element of $B$ over $B'_{n}$, which is given by the identity map $id: C_{n}'\to Z_{n}$. This induces a map $B'_{n}\to B$. We claim that this map is an isomorphism of fibered categories over $\mathfrak{LogSch}$. Given an element of $B$ over $T$ as in (\ref{diag:new-object}), it is equivalent to have a commutative diagram
\[
\xymatrix{
Z_{T} \ar[r] \ar[d] & Z_{n}  \ar[d] \\
C'_{T} \ar[r] \ar[d] \ar@/^1pc/[u]^{s_{T}} & C'_{n} \ar[d] \\
T \ar[r]  & B'_{n} .
}
\]
where the two squares are cartesian squares of log schemes. Note that having $s_{T}$ is equivalent to have a commutative diagram
\[
\xymatrix{
&\underline{f}_{T}^{*}\cM_{Y'} \ar[ld] \ar[rd]& \\
\underline{f}_{T}^{*}\cM_{Y} \ar[rr] && \cM_{C_{n}'}.
}
\]
Since the morphism $\underline{f}_{T}^{*}\cM_{Y'}\to \underline{f}_{T}^{*}\cM_{Y}$ is a surjection, the section $s_{T}$ if exists, is unique. Hence, it is the pull-back of $id: C_{n}'\to Z_{n}$ via $T\to B'_{n}$. This proves that $B = B'_{n}$. 
\end{proof}

\subsection{Minimal objects}
Consider the fibered products (\ref{equ:DF-stack}), (\ref{equ:gen-DF-stack}), or the map (\ref{map:log-embed}). By pulling back the universal maps of each $Y_{i}$, we obtain the universal family of stable Log maps with discrete data $\Gamma$:
 
\begin{equation}\label{diag:univ-family}
\xymatrix{
\cC \ar[r] \ar[d] & Y\\
\cK_{\Gamma}(Y).
}
\end{equation}

This gives the first definition of the minimal stable Log maps:

\begin{definition}\label{defn:exdef-min-map}
A log map $f:C\to Y$ over a log scheme $S$ is called {\em minimal}, if there exists a map $g:\underline{S}\to \underline{\cK}_{\Gamma}(Y)$ of the underlying structure, such that $f$ is obtained by strict pull-back of (\ref{diag:univ-family}) via $g$.
\end{definition}

\begin{proposition}\label{prop:univ-min-map}
Given any stable Log maps $f:C\to Y$ over a log scheme $S$, there exists (up to a unique isomorphism ) a unique minimal stable Log map $f_{m}:C_{m}\to Y$ over $S_{m}$, and a log map $S\to S_{m}$ such that 
\begin{enumerate}
 \item The log map $f$ is given by the pull-back of $f_{m}$ via $S\to S_{m}$.
 \item The underlying map $\underline{S}\to \underline{S}_{m}$ is the identity.
\end{enumerate}
\end{proposition}
\begin{proof}
This follows from the universal property of $\cK_{\Gamma}(Y)$ as a fibered category over $\mathfrak{LogSch}$.
\end{proof}



\section{The combinatorial description of minimality}

\subsection{The case of Deligne-Faltings pairs}\label{ss:DF-marked-graph}

The main result of \cite{Chen} gives more than just Corollary \ref{Cor:stable-maps} for the case $P = \NN$: we have an explicit combinatorial description of the minimal log structure associated to a given log map, in terms of marked graphs. A similar description is possible in general. We first present here the case $P = \NN^m$.

\subsubsection{The graph} Fix a log map $f:C\to Y$ over $S$, such that $\uS$ is a geometric point. Recall that we have a decomposition of $Y$ as in (\ref{equ:decomp-DF-pair}). Denote by $f_{i}$ the following composition 
\begin{equation}\label{diag:decomp-map}
C \stackrel{f}{\rightarrow} Y\to Y_{i}.
\end{equation}
Let $G_{i}$ be the marked graph of the stable Log map $f_{i}$ as in \cite{Chen}. 

We associate to the log map $f$ a {\em $k$-marked graph $G$}, formed out of the dual graph of the curve $C$ with the following extra data:
\begin{enumerate}
\item $k$ partitions of the vertices of $G$ in two types $V(G) =V^{(i)}_0(G)\sqcup V^{(i)}_1(G)$ according to the log maps $f_{i}$ as in \cite{Chen}, where $i=1,\ldots,k$.
\item $k$ integer weights $c^{(i)}_l\geq 0,  i=1,\ldots,m$ on the edges $l\in E(G)$, such that $c^{i}_{l}$ is contact order of $f_{i}$ along $l$.
\item $k$ orientations of edges: whenever an edge $l$ has extremities $v,v'$ and $c^{(i)}_l>0$, we choose one of the orientations $v>_iv'$ or $v'>_iv$ of the edge $l$, which is given by the corresponding orientation in $G_{i}$.
\end{enumerate}

\subsubsection{The monoid}\label{sss:char-monoid}
We introduce a variable $e_{l}$ for each edge $l\in E(G)$, and m variables $e_{v}^{(i)}$, $i=1,2,\cdots,m$ for each vertex $v\in V(G)$. Denote by $h_{l,i}$ the edge equations $e^{(i)}_{v'} = e^{(i)}_v + c^{(i)}_l e_l$ for every edge $l$ with extremities $v\leq_i v'$; and $h_{v,i}$ the vertex equations $e^{(i)}_v=0$, for $v\in V^{(i)}_0(G)$. Consider the monoid
\[
M(G)=\Big\langle\ e_{v},e_{l}\ \Big|\  v\in V_{}(G),\ l\in E(G) \ \Big\rangle\ \bigg/ \big\langle\ h_{l,i}, h_{v,i}\ \big|\  l\in E(G),\ v\in V_{0}^{(i)}(G_{i}) \  \big\rangle
\]
Denote by $T(G)$ the torsion part of $M(G)^{gp}$. Then we have the following composition
\[M(G)\to M(G)^{gp} \to M(G)^{gp}/T(G).\]
Denote by $N(G)$ the image of $M(G)$ in $M(G)^{gp}/T(G)$, and $\ocM(G)$ the saturation of $N(G)$ in $M(G)^{gp}/T(G)$.


\begin{rem}
When $m=1$, the description in \ref{ss:DF-marked-graph} is identical to the case in \cite{Chen}.
\end{rem}

\begin{proposition}\label{prop:DF-monoid}
There is a canonical morphism  $\ocM(G) \to \ocM_S$. A log structure is minimal if and only if this morphism is an isomorphism.
\end{proposition} 
\begin{proof}
Consider the minimal stable Log map $f^{m}:C^{m}\to Y$ over $S^{m}$, and a log map $S\to S^{m}$, which satisfy the conditions in Proposition \ref{prop:univ-min-map}. Denote by $f_{i}$ the stable Log map given by the composition (\ref{diag:decomp-map}). Consider the minimal stable Log map $f_{i}^{min}:C_{i}^{min}\to Y_{i}$ over $S_{i}^{min}$ and the log map $S^{m}\to S_{i}^{min}$ given by Proposition \ref{prop:univ-min-map}. By (\ref{equ:DF-stack}), we have a fiber product of log schemes:
\begin{equation}\label{equ:DF-pt-decomp}
S\cong S_{1}^{min}\times_{\underline{S}}\cdots\times_{\underline{S}}S_{m}^{min}.
\end{equation}
Note that the characteristic monoid of the right hand side of (\ref{equ:DF-pt-decomp}) is given by $\ocM(G)$. Thus, we obtain the map 
\[
\ocM(G) \to \ocM_S.
\]
Assuming the log map $f$ is minimal, this is equivalent to the map $S\to S^{m}$ is an isomorphism, which is equivalent to the map on the level of characteristic $\ocM_{S^{m}}\to \ocM_{S}$ is an isomorphism of monoid. This proves the statement. 
\end{proof}

\subsection{The case of generalized Deligne-Faltings pairs}\label{ss:GDF-marked-graph}
Consider a log map $f:C\to Y$ over $S$, such that $\underline{S}$ is a geometric point.  We now assume that $Y$ is a generalized DF pair with a fixed global presentation $P \to \cM_{Y}$. Thus by (\ref{diag:log-sch-car}), we have a cartesian diagram of log schemes:
\begin{equation}
\xymatrix{
Y \ar[r]^{h} \ar[d]_{h} & Y_{2} \ar[d]^{u_{1}}\\
Y_{1} \ar[r]_{u_{2}} & Y_{3}
}
\end{equation}
where $Y_{i}$ are DF pairs for $i=1,2,3$. Let $f_{i}$ be the following composition 
\begin{equation}
C \stackrel{f}{\rightarrow} Y\to Y_{i}, \ \ \ \mbox{for } i = 1,2,3.
\end{equation}
Then we obtain the marked graph $G_{i}$ for the stable Log map $f_{i}$ as in Section \ref{ss:DF-marked-graph}. Note that by removing all the weights and orientations, the underlying graph $\underline{G}_{i}$ is the dual graph $\underline{G}$ of the underlying curve $C$ over $S$, for all $i$. Denote by $f^{m}_{i}:C_{i}^{m}\to Y_{i}$ over $S_{i}$ the associated minimal log maps of $f_{i}$, and by $f^{m}:C^{m}\to Y$ over $S^{m}$ the associated minimal log maps of $f$. By (\ref{equ:gen-DF-stack}), we obtain a cartesian diagram of log schemes:
\begin{equation}\label{diag:gen-DF-local}
\xymatrix{
S^{m} \ar[d] \ar[r] & S_{2} \ar[d]\\
S_{1} \ar[r] & S_{3}.
}
\end{equation}
Let $\ocM$ be the push-out of the following diagram in the category of toric sharp monoids:
\begin{equation}\label{diag:push-out-base}
\xymatrix{
\ocM(G_{3}) \ar[d] \ar[r] & \ocM(G_{2}) \\
\ocM(G_{1}).  
}     
\end{equation}
 Then we have:

\begin{prop}\label{prop:GDF-monoid}
There is a canonical map $\phi:\ocM\to \ocM_{S}$. The log map $f$ is minimal if and only if $\phi$ is an isomorphism.
\end{prop}
\begin{proof}
The statement can be proved similarly using the argument as for Proposition \ref{prop:DF-monoid}, and the fact that (\ref{diag:gen-DF-local}) is cartesian.
\end{proof}

\begin{remark}
We feel that the expression of $\ocM$ using the push-out diagram (\ref{diag:push-out-base}) is good enough for the purpose of calculations. But it is possible to write down an explicit formation of the base monoids as in Section \ref{sss:char-monoid}. We leave it for interested readers.
\end{remark}

\subsection{Compatibility with Gross-Siebert's construction}
As mentioned in the introduction, another construction of the stack of
stable log maps is given by Gross and Siebert using so called basic
stable log maps when the target is equipped with Zariski log
structures. In fact, when both the Zariski condition and the generalized
Deligne-Faltings or the condition in Proposition \ref{Prop:moreDF}
apply to the log structures on the targets, the notions of basic log
maps and minimal log maps are identical. 

\begin{prop}\label{prop:compare-monoid}
A stable log map is minimal if and only if it is basic in the sense of \cite[Definition 1.17]{GS2}. 
\end{prop}
\begin{proof}

This follows since basic log maps and minimal log maps satisfy the
same universal property: by Proposition \ref{prop:univ-min-map}, any
given log map is uniquely the pullback of 
a unique minimal log map on the same underlying base scheme.  The same
universal property holds for basic log maps by \cite[Proposition 1.20]{GS2}. 
\end{proof}

\begin{rem}
Both minimality and basicness can be defined by putting
constrains on the characteristic monoids of the bases. Thus,
Proposition \ref{prop:compare-monoid} can be proved directly by
comparing the base monoids of both minimal and basic stable log
maps. We provide an argument when the log structure on the target
$Y$ is  Deligne-Faltings. The general situation can be proved
similarly, but with more complicated notation and combinatorics.  

We assume that there is a map $\NN^{k}\to \ocM_{Y}$ locally
lifting to a chart of $\cM_{Y}$. Consider the log scheme $Y_{i}$ as in
Section \ref{ss:DF-marked-graph}. Let $D_{i}$ be the locus with
non-trivial log structure in $Y_{i}$. Consider a stable log map
$f:C\to Y$ over a log scheme $S$. Since both minimality and basicness
can be defined fiber-wise, we might assume that $\uS$ is a geometric
point.   

Consider a generic point $\eta \in \uC$ with associated irreducible
component $v\in V(\uG)$, where $G$ is the marked graph of $f$. In
\cite[Construction 1.15]{GS2}, the factor $P_{\eta}$ can be viewed as
the free monoid generated by the degeneracy $e_{v,i}$ of $v$ for
$i=1,2,\cdots, k$ with the condition that $e_{v,j}=0$ if and only if
the component $v$ does not map into $D_{i}$, or equivalently $v\in
V^{j}_{nd}(G)$. Consider the monoid $\prod_{q\in \uC}\NN$ appearing as
a factor in the expression
\cite[(1.15)]{GS2}, where  $q\in C$ denotes the nodes. It can be viewed as the
free monoid generated by the elements $e_{l}$ for each $l\in
E(G)$. Note that the condition $a_{q}(m)$ in \cite[Construction
  1.15]{GS2} is exactly the the edge condition $h_{l,i}$ in Section
\ref{sss:char-monoid}. Thus, the description of $P_{\eta}$ and
$\prod_{q\in \uC}\NN$ using the elements associated to vertices and
edges induces a natural isomorphism $Q \to \ocM(G)$, where $Q$ is the
monoid defined in \cite[(1.14)]{GS2}. In fact, the object
$(Q,\ocM_{C},\psi,\phi)$ defined in \cite[Construction 1.15]{GS2} is
equivalent to the data of a marked graph $G$. This provides an
explicit derivation of 
Proposition \ref{prop:compare-monoid} in this case.  
\end{rem}

\section{The case of a degeneration}

\subsection{Stable log maps relative to a base}
Consider a family of projective log schemes $\pi: X\to B$, such that $\cM_{X}$ and $\cM_{B}$ are generalized DF log structures. We defined in \cite[Definition 2.1.2]{Chen} that {\em a family of pre-stable log maps over $S$ with target $X/B$} is a commutative diagram of log schemes:
\begin{equation}\label{diag:family-target}
\xymatrix{
C \ar[r]^{f} \ar[d] & X \ar[d] \\
S \ar[r]^{\phi} & B,
}
\end{equation}
such that the family $C\to S$ is a pre-stable log curve. For simplicity, we denote it by $\xi=(C\to S, f, \phi)$, and omit the target $X\to B$, if there is no danger of confusion. The log map $\xi$ is called {\em stable}, if the underlying map $\underline{\xi}$ is stable in the usual sense. 

Consider two pre-stable log maps $\xi_{1}=(C_{1}\to S_{1}, f_{1}, \phi_{1})$ and $\xi_{2}=(C_{2}\to S_{2}, f_{2}, \phi_{2})$. An {\em arrow} $\xi_{1}\to \xi_{2}$ is given by a pair $(g:C_{1}\to C_{2}, h:S_{1}\to S_{2})$ which fits in the following commutative diagram:
\[
\xymatrix{
C_{1} \ar[r]^{g} \ar[d] & C_{2} \ar[r]^{f_{1}} \ar[d] & X \ar[d] \\
S_{1} \ar[r]^{h} & S_{2} \ar[r]^{\phi_{1}} & B,
}
\]
such that the square on the left is a cartesian square of log schemes, and $f_{2}=f_{1}\circ g$ and $\phi_{2}=\phi_{1}\circ h$.

We fix a curve class $\beta$ on the fiber of $\pi$. Define $\cK_{g,n}(X/B,\beta)$ (respectively $\cK_{g,n}^{pre}(X/B,\beta)$) to be the fibered category over $\mathfrak{LogSch}$, parameterizing stable (respectively pre-stable) log maps to $X/B$ with genus $g$, $n$-marked points and curve class $\beta$ on the fiber. 

\subsection{The case when $\pi:X\to B$ is strict}

Note that we have a natural map $\cK_{g,n}(X/B,\beta)\to \cK_{g,n}(\uX/\uB,\beta)$ by removing all log structures from the target. Thus, the stack $\cK_{g,n}(\uX/\uB,\beta)$ is the stack of usual stable maps to $\uX/\uB$ with the canonical log structure from the underlying curves. Note that the stack $\cK_{g,n}(\uX/\uB,\beta)$ is proper over $\uB$. Similarly, we have the natural map $\cK^{pre}_{g,n}(X/B,\beta)\to \cK^{pre}_{g,n}(\uX/\uB,\beta)$, where $\cK^{pre}_{g,n}(\uX/\uB,\beta)$ is the stack parameterizing usual pre-stable maps to $\uX/\uB$ with the canonical log structure from the underlying curves. We first consider the strict case:

\begin{lem}\label{lem:strict-target-base}
Assume that the map $\pi: X\to B$ is strict. Then there is a canonical isomorphism of log stacks 
\[\cK^{pre}_{g,n}(X/B,\beta)\cong \cK^{pre}_{g,n}(\uX/\uB,\beta)\times_{\uB} B.\] 
In particular, by requiring the stability conditions, we have 
\[\cK_{g,n}(X/B,\beta)\cong \cK_{g,n}(\uX/\uB,\beta)\times_{\uB} B.\]
\end{lem}
\begin{proof}
Consider the following commutative diagram
\begin{equation}\label{diag:strict-target-base}
\xymatrix{
\cK^{pre}_{g,n}(X/B,\beta) \ar[rr] \ar[d] && B \ar[d]^{\pi} \\
\cK^{pre}_{g,n}(\uX/\uB,\beta) \ar[rr] && \uB.
}
\end{equation}
We need to prove the above diagram is cartesion. Consider an object 
\[\xi=(p: C\to S, f: C\to X, \phi: S\to B)\in \cK_{g,n}(X/B,\beta).\] 
This is equivalent to have 
\[\xi'=(p:C\to S, \uf: C\to \uX, \underline{\phi}: S\to \uB)\in \cK_{g,n}(\uX/\uB,\beta)\] 
and maps of log structures $f^{\flat}:f^{*}\cM_{X}\to \cM_{C}$ and $\phi^{\flat}:\phi^{*}\cM_{B}\to \cM_{S}$ such that the following diagram of log structures is commutative:
\begin{equation}\label{diag:log-commute}
\xymatrix{
\uf^{*}\circ\pi^{*}\cM_{B} \ar[rr] \ar[d] && \uf^{*}\cM_{X} \ar[d] \\
p^{*}\cM_{S} \ar[rr] && \cM_{C}.
}
\end{equation}
Since the map $\pi$ is strict, the data of $f^{\flat}$, $\phi^{\flat}$, and (\ref{diag:log-commute}) is equivalent to give a map of log schemes $\phi: S\to B$, whose underlying structure is compatible with $\underline{\phi}$ in $\xi'$. This proves that (\ref{diag:strict-target-base}) is cartesian.

The second statement follows since 
$$
\cK_{g,n}(X/B,\beta)\cong\cK_{g,n}(\uX/\uB,\beta)\times_{\cK^{pre}_{g,n}(\uX/\uB,\beta)}\cK^{pre}_{g,n}(X/B,\beta)
$$
\end{proof}

\subsection{Stack parameterizing log sources and targets}

Denote by $\mathfrak{M}_{g,n}$ the stack of pre-stable curves with its canonical log structure. It is proved in \cite{FKato} that the log stack $\mathfrak{M}_{g,n}$ represents the category of all genus $g$, $n$-marked pre-stable log curves over $\mathfrak{LogSch}$. Consider the new stack 
\[\mathfrak{B}=B\times \mathfrak{M}_{g,n}. \]
It represents a fibered category over $\mathfrak{LogSch}$, such that for each log scheme $S$, it associates diagrams of the following form:
\begin{equation}\label{diag:source-target}
\xymatrix{
C \ar[d] \\
S \ar[r]^{\phi} & B,
}
\end{equation}
where $C\to S$ is a log curve. Denote (\ref{diag:source-target}) by $\zeta=(C/S,\phi)$. Given two objects $\zeta_{1}= (C_{1}/S_{1},\phi_{1})$ and $\zeta_{2}=(C_{2}/S_{2},\phi_{2})$, an arrow $\zeta_{2}\to \zeta_{1}$ is given by the following commutative diagram:
\begin{equation}\label{diag:source-target-arrow}
\xymatrix{
C_{1} \ar[r]^{g} \ar[d] & C_{2} \ar[d]  \\
S_{1} \ar[r]^{h} & S_{2} \ar[r]^{\phi_{1}} & B,}
\end{equation}
such that the square is catesian of log schemes, and $\phi_{2}=\phi_{1}\circ h$.

\begin{lem}\label{lem:trivial-fiber-target}
We have the following canonical isomorphism of log stacks:
\[\mathfrak{B} \cong \cK^{pre}_{g,n}(B/B,0).\]
\end{lem}
\begin{proof}
By Lemma \ref{lem:strict-target-base}, we have
\[\cK^{pre}_{g,n}(B/B,0) \cong \cK^{pre}_{g,n}(\uB/\uB,0)\times_{\uB}B\cong \mathfrak{M}_{g,n}\times_{\uB}B\cong \mathfrak{B}.\]
\end{proof}

\subsection{Construction of $\cK_{g,n}(X/B,\beta)$}
We have the following:

\begin{proposition}\label{prop:relate-base-stack}
The fibered categories $\cK_{g,n}(X/B,\beta)$ and $\cK_{g,n}^{pre}(X/B,\beta)$ are represented by algebraic stacks with a natural fs log structures. Furthermore, the underlying stack of $\cK_{g,n}(X/B,\beta)$ is a DM-stack of finite type.
\end{proposition}
\begin{proof}
Consider the stack $\cK_{g,n}^{pre}(B,0)$ as in Remark \ref{rem:all-contact} and the natural map $\cK^{pre}_{g,n}(X,\beta) \to \cK_{g,n}^{pre}(B,0)$ induced by $\pi:X\to B$  as $\pi_*\beta=0$. We have the following commutative diagram:
\begin{equation}\label{diag:relate-base-stack}
\xymatrix{
\cK^{pre}_{g,n}(X/B,\beta) \ar[rr] \ar[d] && \cK^{pre}_{g,n}(X,\beta) \ar[d] \\
\mathfrak{B} \ar[rr] && \cK_{g,n}^{pre}(B,0),
}
\end{equation}
where the right arrow is induced by $\pi$, the left and top arrow is obtained by removing the maps to $X$ and to $B$ in (\ref{diag:family-target}) respectively, and the bottom arrow is obtained by the composition in (\ref{diag:source-target}).

In fact, giving an object $\xi=(C/S,f,\phi)\in \cK^{pre}_{g,n}(X/B,\beta)$ over $S$ is equivalent to giving an object $\zeta=(C/S,\phi)$, and a stable log map $f: C\to X$, which induce the same map $C\to B$. Thus (\ref{diag:relate-base-stack}) is a fibered diagram of fs log stacks. This proves that $\cK^{pre}_{g,n}(X/B,\beta)$ is an algebraic stack with a natural fs log structure. Note that with stability condition, $\cK_{g,n}(X/B,\beta)$ forms an open substack of $\cK^{pre}_{g,n}(X/B,\beta)$, hence is also algebraic.

Note that the image of the $\cK_{g,n}(X/B,\beta)$ in $\mathfrak{B}$ is contained in an open substack of finite type. Therefore, the stack $\cK_{g,n}(X/B,\beta)$ is of finite type. 

Finally, since $\mathfrak{B}\cong \cK_{g,n}^{pre}(B/B,0)$ by Lemma \ref{lem:trivial-fiber-target}, it follows that the bottom arrow of (\ref{diag:relate-base-stack}) is representable. Since the underlying stack of $\cK_{g,n}(X,\beta)$ is a DM-stack, the finiteness of saturation implies that the underlying stack of $\cK_{g,n}(X/B,\beta)$ is also a DM-stack. 


This finishes the proof of the statement.
\end{proof}

\begin{proposition}\label{prop:valuative}
The stack $\cK_{g,n}(X/B,\beta)$ is proper over $\uB$.
\end{proposition}
\begin{proof}
Denote by $\overline{X}=\uX\times_{\uB}B$. Thus we have a canonical map $X\to \overline{X}$, and a strict map $\pi': \overline{X}\to B$. Note also that the log structure on $\overline{X}$ is of generalized Deligne-Faltings. By Proposition \ref{prop:relate-base-stack} and (\ref{diag:relate-base-stack}), we have the following cartesian diagram:
\begin{equation}\label{diag:properness-base}
\xymatrix{
\cK_{g,n}(X/B,\beta) \ar[rr] \ar[d] && \cK_{g,n}(X,\beta) \ar[d] \\
\cK_{g,n}(\overline{X}/B,\beta) \ar[rr] \ar[d] && \cK_{g,n}(\overline{X},\beta) \ar[d] \\
\fB \ar[rr] && \cK_{g,n}^{pre}(B,\pi_{*}\beta).
}
\end{equation}
Note that by Corollary \ref{Cor:stable-maps} both $\cK_{g,n}(X,\beta)$ and $\cK_{g,n}(\overline{X},\beta)$ are representable and finite over $\cK_{g,n}(\uX,\beta)$. It follows from \cite[Lemma 3.12]{L-MB} that the canonical map $\cK_{g,n}(X,\beta)\to\cK_{g,n}(\overline{X},\beta)$ is representable, hence is proper. By Lemma \ref{lem:strict-target-base}, the stack $\cK_{g,n}(\overline{X}/B,\beta)$ is proper over $B$. Thus the statement of the proposition follows.
\end{proof}

\subsection{Minimal objects in the degenerated case}
By (\ref{diag:relate-base-stack}), we have a universal diagram of stable log maps:
\begin{equation}\label{diag:univ-rel-base-map}
\xymatrix{
\cC \ar[rr]  \ar[d] && X \ar[d] \\
\cK_{g,n}(X/B,\beta) \ar[rr] && B.
}
\end{equation}

\begin{definition}\label{defn:min-map-rel-base}
Consider a stable log map $\xi=(C\to S,f,\phi)$ as in (\ref{diag:family-target}). It is called {\em minimal} if there is a map of underlying structures $g:\uS\to \underline{\cK}_{g,n}(X/B,\beta)$, such that $\xi$ is obtained by the strict pull-back of (\ref{diag:univ-rel-base-map}) via $g$.
\end{definition}

\begin{corollary}\label{cor:min-rel-uni}
Given a stable Log map $\xi=(C\to S,f,\phi)$, there exists a minimal stable log map $\xi_{min}=(C_{min} \to S_{min}, f_{min},\phi_{min})$, and a log map $g:S\to S_{min}$ such that 
\begin{enumerate}
 \item $\xi$ is obtained by pull-back $\xi_{min}$ via $g$.
 \item The underlying map $\underline{g}$ is an identity.
\end{enumerate}
Furthermore, the pair $(\xi_{min},g)$ is unique up to a unique isomorphism.
\end{corollary}
\begin{proof}
By Proposition \ref{prop:relate-base-stack}, the log map $\xi$ is equivalent to a log map $S\to \cK_{g,n}(X/B,\beta)$. Define $\xi_{min}$ to be the minimal log map given by the underlying map $\uS\to \underline{\cK}_{g,n}(X/B,\beta)$. This proves the statement.
\end{proof}



\subsection{Compatibility of  minimality}
We show that  minimality in the degeneration case is equivalent to minimality of the map to the total space. First a lemma: 

\begin{lemma}\label{lem:image-B}
Given a stable log map $\xi=(C\to S,f,\phi)\in \cK_{g,n}(X/B,\beta)$, its image in $\cK^{pre}_{g,n}(B,0)$ is a log map with zero contact orders.
\end{lemma}
\begin{proof}
Note that the map $\cK_{g,n}(X/B,\beta)\to \cK^{pre}_{g,n}(B,\pi_{*}\beta)$ factors through $\cB$, which induces stable log maps with only zero contact orders.
\end{proof}

\begin{proposition}\label{prop:min-condition}
Consider $\xi=(C\to S,f,\phi)\in \cK_{g,n}(X/B,\beta)$ a stable Log map over $S$. It is minimal in the sense of Definition \ref{defn:min-map-rel-base}, if and only if the induced stable log map $f:C\to X$ over $S$ in $\cK_{g,n}(X,\beta)$ is minimal in the sense of Definition \ref{defn:exdef-min-map}. 
\end{proposition}
\begin{proof}
In fact, we will prove that the top arrow in (\ref{diag:relate-base-stack}) is strict. This is equivalent to showing that the bottom arrow in (\ref{diag:relate-base-stack}) is strict. This can be checked directly using the construction in the push-out diagram (\ref{diag:push-out-base}).

Consider a diagram 
\[
\xymatrix{
C \ar[d] \\
S \ar[r]^{\phi} & B,
}
\]
induced by a strict map $S\to \cB$. It is enough to consider the case that $\uS$ is a geometric point. Then we obtain an induced log map $h:C\to B$ over $S$. Consider the minimal log map $h':C'\to B$ over $S'$ with the following commutative diagram
\begin{equation}\label{diag:deg-min}
\xymatrix{
C \ar[r]^{j} \ar[d] & C' \ar[r]^{h'} \ar[d]^{p} & B\\
S \ar[r]^{k} & S' \ar@{-->}[ur]_{\underline{\phi}}
}
\end{equation}
such that $h=h'\circ j$, the underlying map $\underline{k}$ is identity, the square is cartesian of fs log schemes, and $\underline{\phi}$ is the underlying map of $\phi$. It is enough to show that the map on the level of characteristic $\bar{k}:\ocM_{S'}\to \ocM_{S}$ is an isomorphism. 

Note that the map $h'$ has only zero contact orders. Otherwise, the composition $h=h'\circ j$ will have non-zero contact orders, which violates Lemma \ref{lem:image-B}. Since all nodes of $C$ is non-distinguished, all the edge equations in Section \ref{sss:char-monoid} is trivial. Hence the construction in Section \ref{ss:GDF-marked-graph} implies a natural splitting $\ocM_{S'}=\ocM_{\uS}^{\uC/\uS}\oplus \ocM$, where $\ocM_{\uS}^{\uC/\uS}$ can be viewed as the submonoid in $\ocM_{S'}$ generated by the elements associated to edges, and $\ocM$ is generated by the element associated to vertices. 

On the other hand, we have $\cM_{S}=\ocM_{\uS}^{\uC/\uS}\oplus \underline{\phi}^{*}\ocM_{B}$. The map $\bar{k}$ induces a map $g: \ocM\to \underline{\phi}^{*}\ocM_{B}$. Using the fact that all contact orders are zero and the construction in (\ref{diag:push-out-base}), we can check that $g$ is an isomorphism. This implies that $\bar{k}$ is an isomorphism.
\end{proof}

\end{document}